\documentclass[preprint,12pt,authoryear]{elsarticle}

\usepackage{graphicx}
\usepackage{epstopdf, epsfig}

\usepackage{color}
\usepackage{amsmath}
\usepackage{amsfonts}
\usepackage{amssymb}
\usepackage{amsbsy}
\usepackage{bm}
\usepackage{graphicx}
\usepackage{extarrows}
\usepackage{tikz}
\usepackage{hyperref}
\usepackage{pgfplots}

\newcommand{\formComma}{\,\text{,}}

\newcommand{\R}{\mathbb{R}}%
\newcommand{\xb}{\mathbf{x}}
\newcommand{\vb}{\mathbf{v}}%
\newcommand{\wb}{\mathbf{w}}%
\newcommand{\eb}{\mathbf{e}}%
\newcommand{\mb}{\mathbf{m}}%
\newcommand{\Xb}{\mathbf{x}}
\newcommand{\surf}{\mathcal{S}}
\newcommand{\gaussianCurvature}{\kappa}
\newcommand{\Grad}{\operatorname{grad}}
\newcommand{\Div}{\operatorname{div}}%
\newcommand{\Rot}{\operatorname{rot}}%
\newcommand{\DivSurf}{\Div_{\surf}}%
\newcommand{\GradSurf}{\Grad_{\surf}}
\newcommand{\RotSurf}{\Rot_{\surf}}
\newcommand{\laplaceBeltrami}{\Delta_{\surf}}
\newcommand{\vecLaplace}{\boldsymbol{\Delta}}
\newcommand{\laplaceDeRham}{\vecLaplace^{\textup{dR}}}
\newcommand{\laplaceRotRot}{\vecLaplace^{\textup{RR}}}
\newcommand{\laplaceGradDiv}{\vecLaplace^{\textup{GD}}}
\newcommand{\LaplaceDeRham}{Laplace-deRham }
\newcommand{\ProjectSurf}{\pi_\surf}
\newcommand{\Tangent}{\mathsf{T}}
\newcommand{\Fs}{\mathcal{T}}

\newcommand{\ie}{i.\,e.}%
\newcommand{\dif}{\textup{d}}
\newcommand{\dS}{\,\dif{\surf}}
\newcommand{\surfNormal}{\boldsymbol{\nu}}
\newcommand{\surfNormalI}{\nu}
\newcommand{\meanCurvature}{\mathcal{H}}
\newcommand{\laplaceDeRhamTilde}{\widehat{\vecLaplace}^{\textup{dR}}}
\newcommand{\lu}{\theta}
\newcommand{\lv}{\varphi}
\newcommand{\extend}[1]{\widehat{#1}}
\newcommand{\extendDomain}[1]{\widetilde{#1}}
\newcommand{\vExt}{\extend{\vb}}
\newcommand{\wExt}{\extend{\wb}}
\newcommand{\triangulation}{{\surf_h}}
\newcommand{\EuBase}[1]{\,\eb^{#1}}
\newcommand{\insertColorbarVertical}[5]{
	\begin{minipage}{0.6cm}
		\begin{tikzpicture}
			\node (colorbar) at (0,0) {\includegraphics[width=\textwidth]{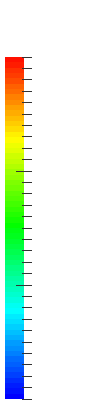}};
			\draw (0.01,1.3) node {\footnotesize #1};
			\draw (-0.15,0.965) node[anchor=west] {\footnotesize #5};
			\draw (-0.15,0.22666667) node[anchor=west] {\footnotesize #4};
			\draw (-0.15,-0.5116667) node[anchor=west] {\footnotesize #3};
			\draw (-0.15,-1.25) node[anchor=west] {\footnotesize #2};
		\end{tikzpicture}
	\end{minipage}
}
\newcommand{\unit}[1]{\mathrm{#1}}
\newtheorem{problem}{Problem}

\begin{document}

\begin{frontmatter}

\title{Solving the incompressible surface Navier-Stokes equation by surface finite elements}

\author[label1]{Sebastian Reuther}
\author[label1,label2,label3]{Axel Voigt}

\address[label1]{Department of Mathematics, TU Dresden, Dresden, Germany}
\address[label2]{Dresden Center for Computational Materials Science (DCMS), Dresden, Germany}
\address[label3]{Center for Systems Biology Dresden (CSBD), Dresden, Germany}

\begin{abstract}
We consider a numerical approach for the incompressible surface Navier-Stokes equation on surfaces with arbitrary genus $g(\surf)$. The approach is based on a reformulation of the equation in Cartesian coordinates of the embedding $\R^3$, penalization of the normal component, a Chorin projection method and discretization in space by surface finite elements for each component. The approach thus requires only standard ingredients which most finite element implementations can offer. We compare computational results with discrete exterior calculus (DEC) simulations on a torus and demonstrate the interplay of the flow field with the topology by showing realizations of the Poincar\'e-Hopf theorem on $n$-tori.
\end{abstract}


\end{frontmatter}

\section{Introduction}
\label{sec1}

We consider a compact smooth Riemannian surface $\surf$ without boundary and an incompressible surface Navier-Stokes equation
\begin{align}
	\label{eq1}
	\partial_{t}\vb + \nabla_{\vb} \vb &= - \GradSurf p + \frac{1}{\text{Re}} \left(- \laplaceDeRham \vb + 2 \gaussianCurvature \vb \right) \\
	\label{eq2}
	\DivSurf \vb &= 0  
\end{align}
in $\surf\times\left( 0,\infty \right)$ with initial condition \( \vb \left( \xb, t=0  \right) = \vb_{0}(\xb) \in \Tangent_{\xb}\surf \). Thereby $\vb(t) \in \Tangent
\surf$ denotes the tangential surface velocity, $p(\xb,t) \in \R$ the surface pressure, $\text{Re}$ the surface Reynolds number, $\gaussianCurvature$ the Gaussian curvature,
$\Tangent_{\xb}\surf $ the tangent space on $ \xb \in \surf$, $\Tangent \surf = \cup_{\xb \in \surf} \Tangent_{\xb} \surf$ the tangent bundle and $\nabla_{\vb}, \DivSurf$ as well as
$\laplaceDeRham$ the covariant directional derivative, surface divergence and surface \LaplaceDeRham operator, respectively. The surface Navier-Stokes equation results from conservation of mass and (tangential) linear momentum. Alternatively, eqs. \eqref{eq1} and \eqref{eq2} can also be derived from the Rayleigh dissipation potential \cite{Doerriesetal_PRE_1996} or as a thin-film limit of the three-dimensional incompressible Navier-Stokes equation \cite{Miura_arXiv_2017}. 

The incompressible surface Navier-Stokes equation is related to the Bouss\-inesq-Scriven constitutive law for the surface viscosity in two-phase flow problems \cite{Scriven_CES_1960,Secombetal_QJMAM_1982,Botheetal_JFM_2010} and to fluidic biomembranes \cite{Huetal_PRE_2007,ArroyoDeSimone_PRE_2009,Barrettetal_PRE_2015,Reutheretal_JCP_2016}. 
Further applications can be found in computer graphics, e.g. \cite{Elcottetal_2007,Mullenetal_2009,Vaxmanetal_2016}, and geophysics, e.g. \cite{Padberg-Gehleetal_TMDAVIV_2017,sasakietal_JFM_2015}.
The equation is also studied as a mathematical problem of its own interest, see e.g. \cite{EbinMarsden_AM_1970,MitreaTaylor_MA_2001}.

While a huge literature exists for the numerical treatment of the two-dimensional incompressible Navier-Stokes equation in flat space, results for its surface counterpart eqs. \eqref{eq1} and \eqref{eq2} are rare. In \cite{Nitschkeetal_JFM_2012,Reutheretal_MMS_2015} a surface vorticity-stream function formulation is introduced. However, this approach cannot deal with harmonic vector fields and is therefore only applicable on surfaces with genus \( g(\surf) = 0 \). The only direct numerical approach for eqs. \eqref{eq1} and \eqref{eq2}, which is also desirable for surfaces with genus $g(\surf) \neq 0$, was proposed in \cite{Nitschkeetal_book_2017} and uses discrete exterior calculus (DEC). 

The purpose of this paper is to introduce a surface finite element discretization with only standard ingredients. This is achieved by extending the variational space from vectors in $\Tangent \surf$ to vectors in $\R^3$ and penalizing the normal component. This allows to split the vector-valued problem into a set of coupled scalar-valued problems for each component for which standard surface finite elements, see the review \cite{Dziuketal_AN_2013}, can be used. Similar approaches have already been independently used for other vector-valued problems, see \cite{Hansboetal_arXiv_2016} for a surface vector Laplacian, \cite{Nestleretal_JNS_2017} for a surface Frank-Oseen problem and \cite{Jankuhnetal_arXiv_2017} for a surface Stokes problem.

The paper is organized as follows. In Section \ref{sec2} we introduce the necessary notation, reformulate the problem in Cartesian coordinates of the embedding $\R^3$ and introduce the penalization of the normal component. We further modify the equation by rotating the velocity field, which reduces the complexity of the equation. In Section \ref{sec3} we describe the numerical approach. For the resulting equations we propose a Chorin projection approach and a discretization in space by standard piecewise linear Lagrange surface finite elements. We demonstrate the reduction of computational time due to the introduced rotation and validate our approach against a DEC solution on a torus with harmonic vector fields, see \cite{Nitschkeetal_book_2017}. In Section \ref{sec4} results are shown and analyzed on $n$-tori and conclusions are drawn in Section \ref{sec5}.

\section{Model formulation}
\label{sec2}

We follow the same notation as introduced in \cite{Nestleretal_JNS_2017} and parametrize the surface $\surf\subset\R^3$ by the local coordinates $ \lu, \lv$, \ie, 
\begin{equation*} 
	\Xb:\R^{2}\supset U \rightarrow \R^{3};\ \left( \lu, \lv \right)\mapsto \Xb \left( \lu, \lv \right)\, .
\end{equation*}
Thus, the embedded \( \R^{3} \) representation of the surface is given by \( \surf = \Xb(U) \). The unit
outer normal of $\surf$ at point $\Xb$ is denoted by $\surfNormal(\Xb)$. We denote by \( \left\{ \partial_{\lu}\Xb,\partial_{\lv}\Xb  \right\}\) the canonical basis to describe the (tangential) velocity
\( \vb(\Xb)\in \Tangent_{\Xb}\surf \), \ie, \( \vb = v^{\lu} \partial_{\lu}\Xb + v^{\lv}\partial_{\lv}\Xb \) at a point \( \Xb\in\surf \).
In a (tubular) neighborhood $\Omega_\delta$ of $\surf$, defined by $\Omega_\delta := \{ \extendDomain{\Xb}\in\R^3\,:\,d_\surf(\extendDomain{\Xb}) < \frac{1}{2}\delta \}$,
with a signed-distance function $d_\surf(\extendDomain{\Xb})$ a coordinate projection $\Xb\in\surf$ of $\extendDomain{\Xb}\in\R^3$ is introduced, such that $\extendDomain{\Xb} = \Xb + d_\surf(\extendDomain{\Xb})\surfNormal(\Xb)$. For $\delta$ sufficiently small (depending on the local curvature of the surface)
this projection is injective, see \cite{Dziuketal_AN_2013}. For a given $\extendDomain{\Xb}\in\Omega_\delta$ the coordinate
projection of $\extendDomain{\Xb}$ will also be called gluing map, denoted by $\pi:\Omega_\delta\to\surf,\,\extendDomain{\Xb}\mapsto\Xb$. The pressure $p:\surf\to\R$ and the velocity $\vb:\surf\to\Tangent\surf$ can be smoothly extended
in the neighborhood $\Omega_\delta$ of $\surf$ by utilizing the coordinate projection, \ie,
extended fields $\tilde{p}:\Omega_\delta\to\R$ and $\tilde{\vb}:\Omega_\delta\to\R^3$ are defined by
\begin{equation}
	\label{eq:smooth_extension}
	\tilde{p}(\extendDomain{\Xb}) := p(\Xb)\quad\text{ and }\quad\tilde{\vb}(\extendDomain{\Xb}) := \vb(\Xb)\formComma
\end{equation}
respectively, for $\extendDomain{\Xb}\in\Omega_\delta$ and $\Xb$ the corresponding coordinate projection. To embed the \( \R^{3} \) vector space structure to the tangential bundle of the surface we use the pointwise defined normal projection
\begin{align*} 
	\begin{aligned}
		\ProjectSurf(\xb):  \Tangent_\xb\R^3 \cong \R^{3} &\rightarrow \Tangent_\xb\surf;\\
		\vExt(\xb) &\mapsto \vExt(\xb) - \surfNormal(\xb)(\surfNormal(\xb)\cdot \vExt(\xb)) = \vb(\xb)
	\end{aligned}
\end{align*}
for all \( \xb \in \surf \), which maps the \( \R^{3} \) velocity
\( \vExt = v_{x}\EuBase{x} + v_{y}\EuBase{y} + v_{z}\EuBase{z} \in \R^3  \), not necessarily tangential to the surface,
to the tangential velocity \( \vb \in \Tangent_\xb\surf \). We drop the argument $\xb$ when applied to velocity fields living on $\surf$.
With these notations we have the following correspondence of the different representations of first order differential operators on surfaces:
\vspace{0.01\textwidth}
\begin{center}
	\begin{tabular}{lllll}
		$\Tangent\surf$ & $\GradSurf p$          & $\RotSurf p$                  &  $\DivSurf \vb$ & $\RotSurf \vb$ \\ \hline
		$\R^3$          &  $\ProjectSurf\nabla p$ &  $\surfNormal \times \nabla p$ & $\nabla \cdot \vExt - \surfNormal \cdot (\nabla \vExt \cdot \surfNormal) - \meanCurvature (\vExt \cdot \surfNormal)$ & $\left( \nabla\times\vExt \right)\cdot\surfNormal$
	\end{tabular}
\end{center}
\vspace{0.01\textwidth}
Thereby $\meanCurvature$ denotes the mean curvature. 
We further define $\DivSurf \vExt = \nabla \cdot \vExt - \surfNormal \cdot (\nabla \vExt \cdot \surfNormal)$, $\RotSurf \vExt = -\DivSurf (\surfNormal \times \vExt)$ and . Using the definition of  \cite{Marsden1988} the surface \LaplaceDeRham operator is defined as $\laplaceDeRham \vb = -\left(\laplaceRotRot  +  \laplaceGradDiv \right) \vb$ with $\laplaceRotRot \vb = \RotSurf\RotSurf\vb$ and $\laplaceGradDiv \vb = \GradSurf\DivSurf\vb$. As shown in \cite{Nestleretal_JNS_2017} it holds 
\begin{align*}
	\laplaceDeRham \vb \approx \laplaceDeRhamTilde \vExt \quad \text{with} \quad \laplaceDeRhamTilde \vExt = - (\RotSurf\RotSurf\vExt + \GradSurf\DivSurf\vExt)
\end{align*} 
if the normal component $(\vExt \cdot \surfNormal)$ is penalized by the additional term $\alpha (\surfNormal \cdot \vExt) \surfNormal$. First order convergence in the penalty parameter $\alpha$ was numerically shown for this approximation in \cite{Nestleretal_JNS_2017}. Due to the incompressibility we here have $\laplaceDeRham \vb \approx -\RotSurf\RotSurf\vExt$. With $\nabla_{\vb} \vb = \frac{1}{2} \GradSurf (\vExt \cdot \vExt) + \RotSurf \vExt \surfNormal \times \vExt$ we thus obtain the approximation of the surface incompressible Navier-Stokes equation in Cartesian coordinates which ensures the velocity to be tangential only weakly through the added penalty term
\begin{align}
	\label{eq3}
	\partial_{t}\vExt + \RotSurf \vExt \surfNormal \times \vExt &= - \GradSurf \tilde{p} + \frac{1}{\text{Re}} \left(- \RotSurf\RotSurf \vExt + 2 \gaussianCurvature \vExt \right) - \alpha (\vExt \cdot \surfNormal) \surfNormal \\
	\label{eq4}
	\DivSurf \vExt &= 0  
\end{align}
with $\tilde{p} = p + \frac{1}{2} \vExt \cdot \vExt$. Eqs. \eqref{eq3} and \eqref{eq4} can now be solved for each component $v_x$, $v_y$, $v_z$ and $\tilde{p}$ using standard approaches for scalar-valued problems on surfaces, such as the surface finite element method \cite{DziukElliott_JCM_2007,DziukElliott_IMAJNA_2007,Dziuketal_AN_2013}, level set approaches \cite{Bertalmioetal_JCP_2001,Greeretal_JCP_2006,Stoeckeretal_JIS_2008,dziuketal_IFB_2008} or diffuse interface approximations \cite{Raetzetal_CMS_2006}. However, the $\RotSurf\RotSurf \vExt$ term leads to a heavy workload in terms of implementation and assembly time, as 36 second order operators, 72 first order operators and 36 zero order operators have to be considered. This effort can drastically be reduced by rotating the velocity field in the tangent plane. Instead of $\vExt$ we consider $\wExt = \surfNormal \times \vExt$ as unknown. Applying $\surfNormal \times$ to eq. \eqref{eq3} we thus obtain
 \begin{align}
	\label{eq5}
	\partial_{t}\wExt + \DivSurf \wExt \surfNormal \times \wExt &= - \RotSurf \tilde{p} + \frac{1}{\text{Re}} \left( \GradSurf\DivSurf \wExt + 2 \gaussianCurvature \wExt \right) - \alpha (\wExt \cdot \surfNormal) \surfNormal \\
	\label{eq6}
	\RotSurf \wExt &= 0  
\end{align}
where we have used the identities $\RotSurf \vExt = - \DivSurf \wExt$, $\DivSurf \vExt = \RotSurf \wExt$, $\vExt = - \surfNormal \times \wExt$ and $\surfNormal \times ( \surfNormal \times \vExt) = - \vExt$. The $\GradSurf\DivSurf \wExt$ term now contains only 9 second order terms and the remaining terms are of similar complexity as in eqs. \eqref{eq3} and \eqref{eq4}.

\section{Discretization}
\label{sec3}

\subsection{Time discretization}

Let $0 < t_0 < t_1 < \ldots$ be a sequence of discrete times with time step width $\tau_n := t_{n+1}-t_n$ in the $n$-th iteration. The fields $\vExt^n(\xb)e.g.uiv \vExt(\xb,t_n)$, $\wExt^n(\xb)e.g.uiv \wExt(\xb,t_n)$ and $\tilde{p}^n(\xb)e.g.uiv \tilde{p}(\xb,t_n)$ correspond to the time-discrete functions at $t_n$. Applying a Chorin projection method \cite{Chorin_MC_1968} to eqs. \eqref{eq3} and \eqref{eq4} with a semi-implicit Euler time scheme results in time discrete systems of equations as follows:

\begin{problem}\label{prob:1}
	Let $\vExt^0\in C(\surf;\,\R^3)$ be a given initial velocity field with $\vExt^0 = \vb^0$. For $n=0,1,2,\ldots$ find 
	\begin{enumerate}
		\item $\vExt^*$ such that
			\begin{align*}
				\!\!\!\!\!\frac{1}{\tau_n} (\vExt^* - \vExt^n) &= - \RotSurf \vExt^* \surfNormal \times \vExt^n + \frac{1}{\text{Re}} \left(- \RotSurf\RotSurf \vExt^* + 2 \gaussianCurvature \vExt^* \right) - \alpha (\vExt^* \cdot \surfNormal) \surfNormal 
			\end{align*}
		\item $\tilde{p}^{n+1}$ such that
			\begin{align*}
				\tau_n \laplaceBeltrami \tilde{p}^{n+1} &= \DivSurf \vExt^*
			\end{align*}
		\item $\vExt^{n+1}$ such that
			\begin{align*}
			\vExt^{n+1} &= \vExt^* - \tau_n \GradSurf \tilde{p}^{n+1},
			\end{align*}
	\end{enumerate}
	with $\laplaceBeltrami$ the Laplace-Beltrami operator. 
\end{problem}

The corresponding scheme for eqs. \eqref{eq5} and \eqref{eq6} follows by defining $\wExt^* = \surfNormal \times \vExt^*$ and applying $\surfNormal \times$ to the equation in the first step. We thus obtain:

\begin{problem}\label{prob:2}
	Let $\vExt^0 \in C(\surf;\,\R^3)$ be a given initial velocity field with $\vExt^0 = \vb^0$. Compute $\wExt^0 = \surfNormal \times \vExt^0$. For $n=0,1,2,\ldots$ find 
	\begin{enumerate}
		\item $\wExt^*$ such that
			\begin{align*}
				\!\!\!\!\!\!\!\!\frac{1}{\tau_n} (\wExt^* - \wExt^n) &= \DivSurf \wExt^* \surfNormal \times \wExt^n + \frac{1}{\text{Re}} \left( \GradSurf\DivSurf \wExt^* + 2 \gaussianCurvature \wExt^* \right) - \alpha (\wExt^* \cdot \surfNormal) \surfNormal 
			\end{align*}
		\item $\tilde{p}^{n+1}$ such that
			\begin{align*}
				\tau_n \laplaceBeltrami \tilde{p}^{n+1} &= \RotSurf \wExt^*
			\end{align*}
		\item $\wExt^{n+1}$ such that
			\begin{align*}
			\wExt^{n+1} &= \wExt^* - \tau_n \RotSurf \tilde{p}^{n+1}
			\end{align*}
		\item $\vExt^{n+1} = - \surfNormal \times \wExt^{n+1}$.
	\end{enumerate}
\end{problem}

For simplicity we consider only a Taylor-$0$ linearization of the nonlinear term in both problems.

\subsection{Space discretization}

For the discretization in space we apply the surface finite element method for scalar-valued problems \cite{Dziuketal_AN_2013} for each component.
Therefore, the surface $\surf$ is discretized by a conforming triangulation $\triangulation$, given as the union of simplices, \ie, $\triangulation := \bigcup_{T \in \Fs} T$.
We use globally continuous, piecewise linear Lagrange surface finite elements
\[
\mathbb{V}_h(\triangulation) = \left\{ v_h \in C^0(\triangulation) \,:\, v_h|_T \in \mathbb{P}^1, \, \forall \, T \in \Fs \right\}
\]
as trial and test space for all components $\extend{v}_{i}$ of $\vExt$ as well as $\extend{w}_{i}$ of $\wExt$ and $\tilde{p}$ with $\Fs$ the set of triangular faces.

The resulting fully discrete problem for \autoref{prob:1} reads: For $n=0,1,2,\ldots$ find $\extend{v}_{i}^*$, $\tilde{p}^{n+1}\in\mathbb{V}_h(\triangulation)$ s.t. $\forall\,\extend{u}_i, \extend{q} \in \mathbb{V}_h(\triangulation)$
\begin{align}
	\!\!\!\!\!\frac{1}{\tau_n} \int_{\surf_h} \extend{v}_{i}^*\extend{u}_i \dS + \int_{\surf_h} \RotSurf \vExt^* \left(\surfNormal \times \vExt^n\right)_i \extend{u}_i \dS + \alpha \int_{\surf_h} \surfNormal &\cdot\vExt^* \surfNormalI_{i} \extend{u}_i \dS \nonumber \\
	- \frac{1}{\text{Re}} \int_{\surf_h} \RotSurf \vExt^* \RotSurf \left(\extend{u}_i \eb^i\right)\dS - 2 \int_{\surf_h}\gaussianCurvature \extend{v}_i^*\extend{u}_i \dS 
	&= \frac{1}{\tau_n} \int_{\surf_h} \extend{v}_{i}^{n}\extend{u}_i \dS \\
	\tau_n \int_{\surf_h} \GradSurf \tilde{p}^{n+1} \cdot \GradSurf \extend{q} \dS + \int_{\surf_h} \vExt^* \cdot \GradSurf \extend{q} \dS &= 0
\end{align}
for $i = x,y,z$, from which $\vExt^{n+1}$ can be computed according to step 3 in \autoref{prob:1}. 

The resulting fully discrete problem for \autoref{prob:2} reads: For $n=0,1,2,\ldots$ find $\extend{w}_{i}^*$, $\tilde{p}^{n+1}\in\mathbb{V}_h(\triangulation)$ s.t. $\forall\,\extend{u}_i, \extend{q} \in \mathbb{V}_h(\triangulation)$

\begin{align}
	\!\!\!\!\!\frac{1}{\tau_n} \int_{\surf_h} \extend{w}_{i}^*\extend{u}_i \dS - \int_{\surf_h} \DivSurf \wExt^* (\surfNormal \times \wExt^n)_i \extend{u}_i \dS + \alpha \int_{\surf_h} \surfNormal &\cdot\wExt^* \surfNormalI_{i} \extend{u}_i \dS \nonumber \\
	+ \frac{1}{\text{Re}} \int_{\surf_h} \DivSurf \wExt^* \left(\GradSurf \extend{u}_i\, \right)_i \dS - 2 \int_{\surf_h}\gaussianCurvature \extend{w}_{i}^*\extend{u}_i \dS 
	&= \frac{1}{\tau_n} \int_{\surf_h} \extend{w}_{i}^{n}\extend{u}_i\dS \\
	\!\!\!\!\tau_n \int_{\surf_h} \GradSurf \tilde{p}^{n+1} \cdot \GradSurf \extend{q} \dS + \int_{\surf_h} \surfNormal \times \wExt^* \cdot \GradSurf \extend{q} \dS &= 0
\end{align}
for $i = x,y,z$, from which $\wExt^{n+1}$ and $\vExt^{n+1}$ can be computed according to step 3 and 4 in \autoref{prob:2}. 

To assemble and solve the resulting system we use the FEM-toolbox AMDiS \cite{VeyVoigt_CVS_2007,Witkowskietal_ACM_2015} with domain decomposition on $16$ processors. As linear solver we have used a BiCGStab($l$) method with $l = 2$ and a Jacobi preconditioner with ILU($0$) local solver on each partition.

\subsection{Comparison and validation}

\begin{figure}
	\begin{center}
		\begin{minipage}{0.57\textwidth}
			\begin{center}
				\begin{tabular}{rrrrrr}
					\multicolumn{1}{c}{number} & \multicolumn{1}{c}{time $t_{\vExt}$} & \multicolumn{1}{c}{time $t_{\wExt}$} & \multicolumn{1}{c}{time ratio} \\
					\multicolumn{1}{c}{of DOFs} & \multicolumn{1}{c}{(in $\unit{ms}$)} & \multicolumn{1}{c}{(in $\unit{ms}$)} & \multicolumn{1}{c}{$t_{\vExt} / t_{\wExt}$} \\
					\hline
					$4614$  &  $1123.39$ &  $15.56$ & $72.2$ \\
					$9222$  &  $2199.54$ &  $31.62$ & $69.6$ \\
					$18438$ &  $4368.73$ &  $69.05$ & $63.3$ \\
					$36870$ &  $8817.76$ & $156.68$ & $56.3$ \\
					$73734$ & $17920.00$ & $326.54$ & $54.9$ \\
				\end{tabular}
			\end{center}
		\end{minipage} \hfill
		\begin{minipage}{0.42\textwidth}
			\begin{center}
%
%
\begin{tikzpicture}

\begin{axis}[%
width=0.8\textwidth,
height=0.6\textwidth,
at={(0,0)},
scale only axis,
separate axis lines,
every outer x axis line/.append style={black},
every x tick label/.append style={font=\color{black}},
xmin=0,
xmax=80,
xtick={20,40,60,80},
xlabel={$\mbox{number of DOFs (in thousands)}$},
every outer y axis line/.append style={black},
every y tick label/.append style={font=\color{black}},
ymin=0,
ymax=20,
ytick={0,5,10,15,20},
ylabel={$\mbox{time }t\mbox{ (in s)}$},
axis background/.style={fill=white},
x label style={font=\footnotesize,at={(axis description cs:0.5,0.05)}},
y label style={font=\footnotesize,at={(axis description cs:0.1,0.5)}},
yticklabel style = {font=\footnotesize},
xticklabel style = {font=\footnotesize},
legend style={font=\footnotesize,at={(0.02,0.98)},anchor=north west,legend cell align=left,align=left,draw=black}
]
\addplot [color=blue,solid,line width=1.0pt,mark=o,mark options={solid}]
  table[row sep=crcr]{%
73.734	17.92\\
36.87	8.81776\\
18.438	4.36873\\
9.222	2.19954\\
4.614	1.12339\\
2.31	0.543686\\
1.158	0.270554\\
0.582	0.138096\\
0.294	0.0697403\\
0.15	0.0342855\\
0.078	0.0230002\\
0.042	0.0091605\\
0.024	0.0049905\\
};
\addlegendentry{$\mbox{rot}_{\mathcal{S}}\mbox{rot}_{\mathcal{S}} \widehat{\mathbf{v}}$};

\addplot [color=red,dashed,line width=1.0pt,mark=o,mark options={solid}]
  table[row sep=crcr]{%
73.734	0.326541\\
36.87	0.156677\\
18.438	0.0690535\\
9.222	0.0316153\\
4.614	0.0155631\\
2.31	0.00843\\
1.158	0.0049032\\
0.582	0.00247\\
0.294	0.0013735\\
0.15	0.0011907\\
0.078	0.0006026\\
0.042	0.000622\\
0.024	0.00048\\
};
\addlegendentry{$\mbox{grad}_{\mathcal{S}}\mbox{div}_{\mathcal{S}} \widehat{\mathbf{w}}$};

\end{axis}
\end{tikzpicture}%
			\end{center}
		\end{minipage} 
	\end{center}
	\caption{Assembly times $t_{\vExt}$ and $t_{\wExt}$ for the two second order operators $\RotSurf\RotSurf \vExt$ and $\GradSurf\DivSurf \wExt$ as a function of the number of DOFs.}
	\label{fig1}
\end{figure}
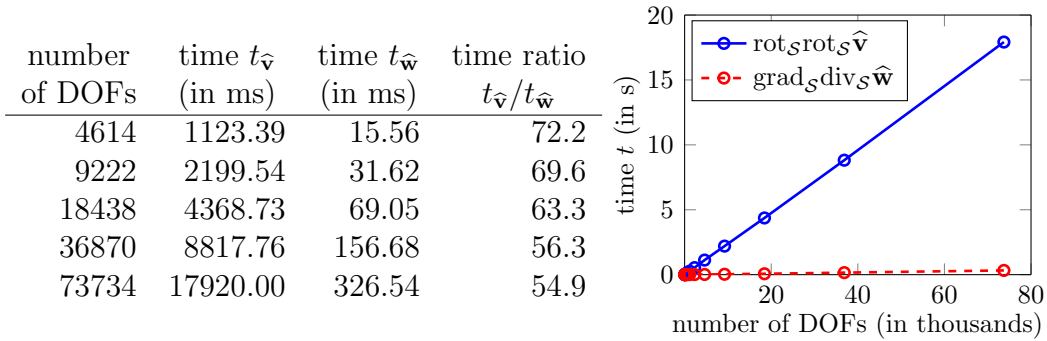

Both approaches lead to the same results. However, the computational cost for \autoref{prob:2} is drastically reduced. To quantify this reduction we compare the assembly time for the second order operators in \autoref{prob:1} and \autoref{prob:2}. We consider a sphere as computational domain $\surf = \mathbb{S}^2$ and vary the triangulation $\Fs$. Figure \ref{fig1} shows the assembly time as a function of degrees of freedom (DOFs). The time is the mean value of multiply runs of the assembly routine. The results indicate a reduction by a factor of approximately 50.

We now compare the solution of \autoref{prob:2} with an example considered in \cite{Nitschkeetal_book_2017} using DEC. It considers a nontrivial solution with $\DivSurf \vb =0$ and $\RotSurf \vb = 0$. Such harmonic vector fields can exist on surfaces with $g(\surf) \neq 0$. We consider a torus which has genus $g(\surf) =1$. A torus can be described by the levelset function $T(\xb) = (x^2+y^2+z^2 + R^2 - r^2)^2 - 4 R^2 (x^2 + z^2)$ with $\xb = (x,y,z) \in \R^3$, major radius $R$ and minor radius $r$. We here use $R = 2$ and $r=0.5$. Let $\phi$ and $\theta$ denote the standard parametrization angles on the torus. Then, the two basis vectors can be written as $\partial_\phi\xb$ as well as $\partial_\theta\xb$ and read in Cartesian coordinates $\partial_\phi\xb = (-z, 0, x)$ as well as $\partial_\theta\xb = (-\frac{xy}{\sqrt{x^2+z^2}},  \sqrt{x^2+z^2} - 2, -\frac{yz}{\sqrt{x^2+z^2}})$. There are two (linear independent) harmonic vector fields on the torus,
\begin{align*}
	\vb_\phi^{harm} = \frac{1}{4\left(x^2+z^2\right)} \partial_\phi\xb \qquad &\mbox{and} \qquad
	\vb_\theta^{harm} = \frac{1}{2\sqrt{x^2+z^2}} \partial_\theta\xb.
\end{align*}
The example considers the mean of the two harmonic vector fields as initial condition $\vb_0(\xb) = \frac{1}{2}( \vb_\phi^{harm} + \vb_\theta^{harm} )$ and shows the evolution towards a Killing vector field which is proportional to the basis vector $\partial_\phi\xb$.
The surface Reynolds number is $\text{Re}=10$. Figure \ref{fig2} shows the results obtained with the fully discrete scheme of \autoref{prob:2} with time step width $\tau_n = 0.1$ and penalization parameter $\alpha = 3000$ on the same mesh as considered in \cite{Nitschkeetal_book_2017}. For the Gaussian curvature $\gaussianCurvature$ we use the analytic formula. 

\begin{figure}
	\begin{center}
		\begin{minipage}{0.17\textwidth}
			\includegraphics[width=\textwidth]{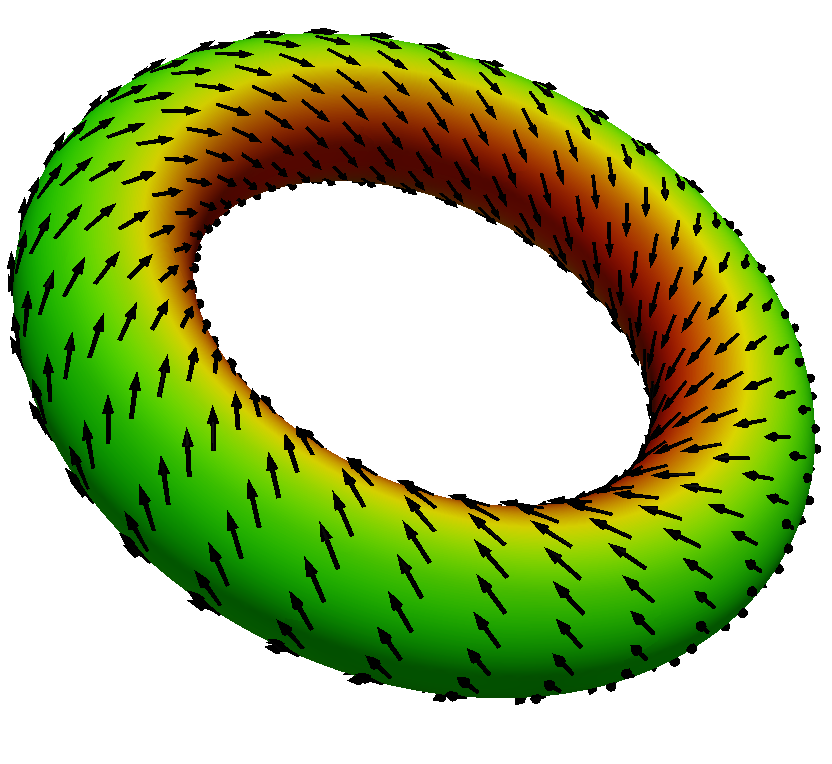}
		\end{minipage}
		\begin{minipage}{0.17\textwidth}
			\includegraphics[width=\textwidth]{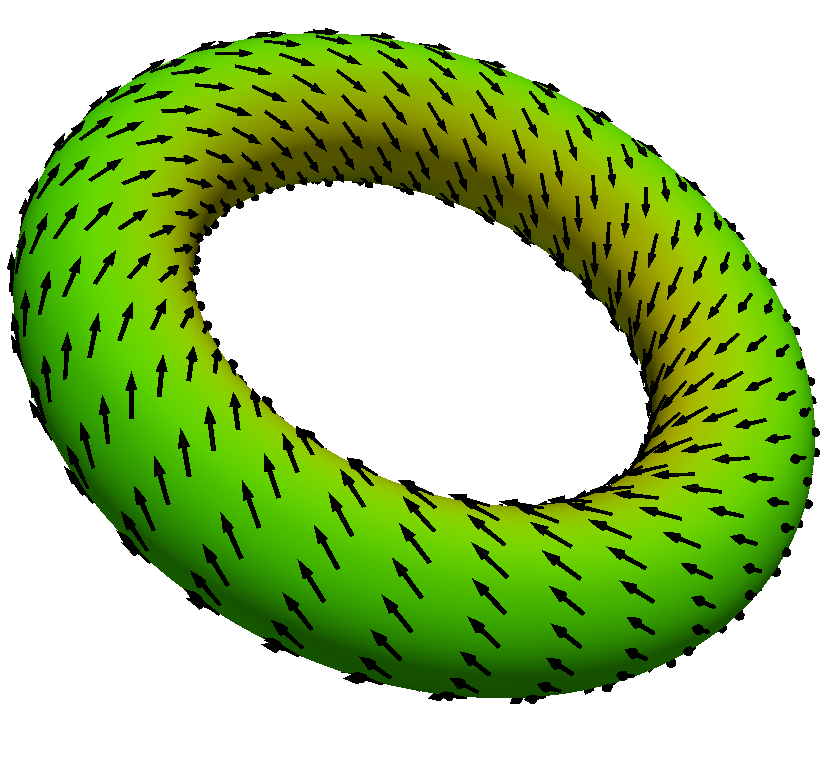}
		\end{minipage}
		\begin{minipage}{0.17\textwidth}
			\includegraphics[width=\textwidth]{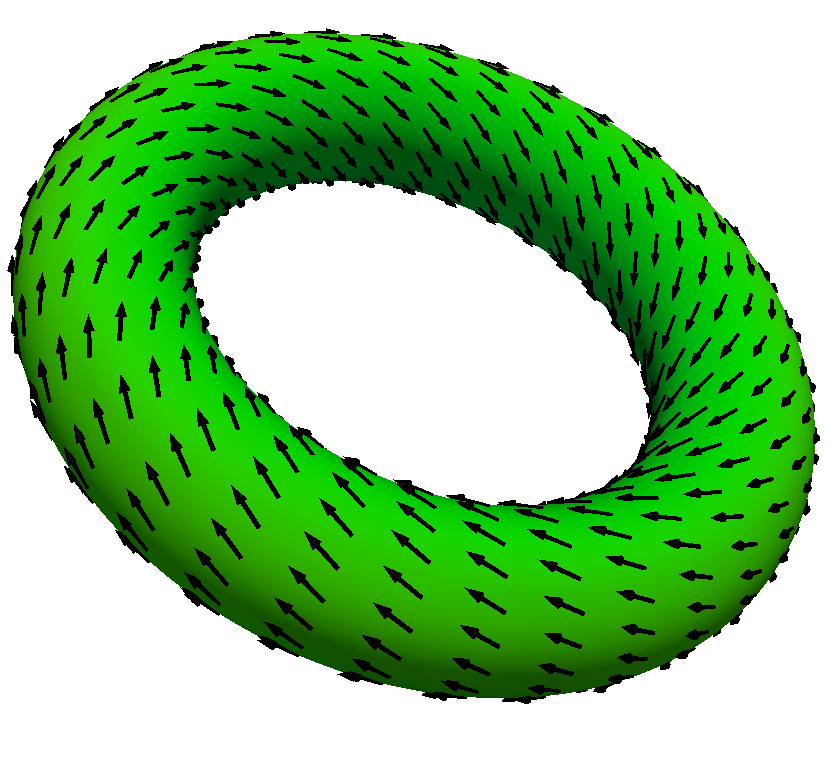}
		\end{minipage}
		\begin{minipage}{0.17\textwidth}
			\includegraphics[width=\textwidth]{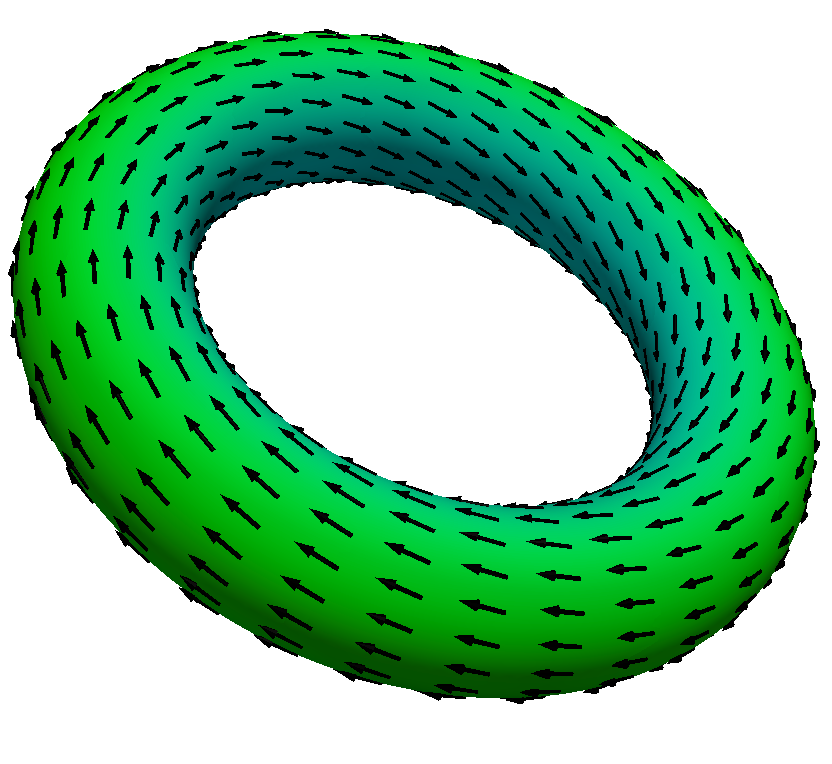}
		\end{minipage}
		\begin{minipage}{0.17\textwidth}
			\includegraphics[width=\textwidth]{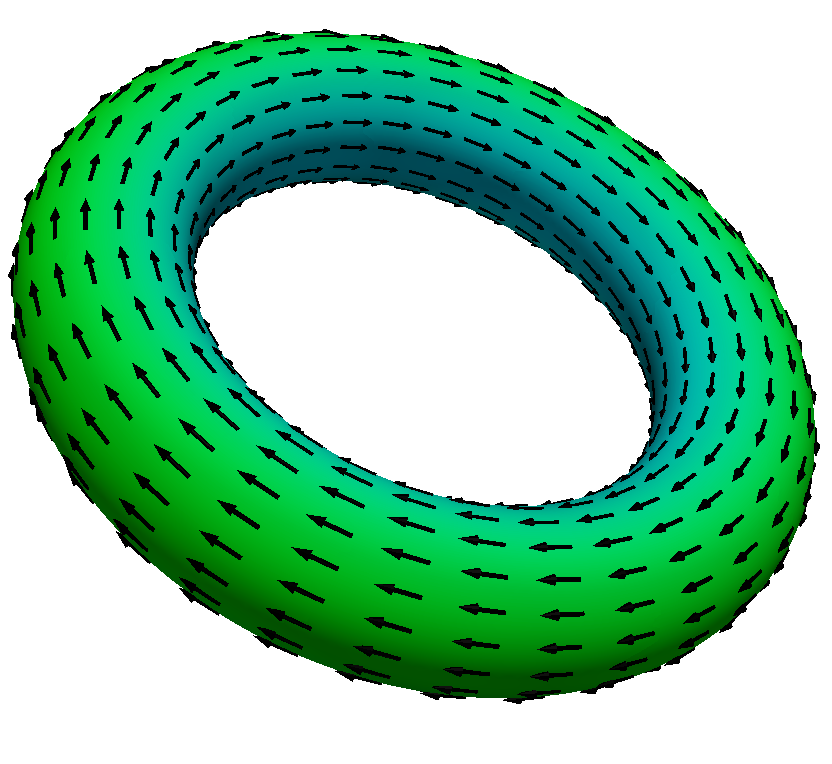}
		\end{minipage}
		\insertColorbarVertical{$\|\vExt\|$}{$0.02$}{}{}{$0.12$}
	\end{center}
	\caption{Numerical solution of $\vExt = - \surfNormal \times \wExt$ at $t = 0$, $2$, $10$, $30$ and $60$ (left to right). The color indicates the absolute value of the velocity $\vExt$. The arrows are rescaled for better visualization.}
	\label{fig2}
\end{figure}
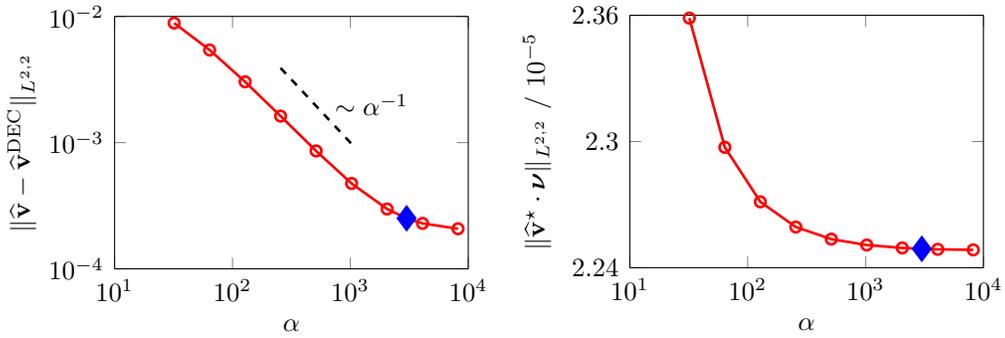
\begin{figure}
	\begin{center}
		\begin{minipage}{0.49\textwidth}
%
%
\begin{tikzpicture}

\begin{axis}[%
width=0.7\textwidth,
height=0.5\textwidth,
at={(0,0)},
scale only axis,
separate axis lines,
every outer x axis line/.append style={black},
every x tick label/.append style={font=\color{black}},
xmode=log,
xmin=10,
xmax=10000,
xtick={   10,   100,  1000, 10000},
xminorticks=false,
xlabel={$\alpha$},
every outer y axis line/.append style={black},
every y tick label/.append style={font=\color{black}},
ymode=log,
ymin=0.0001,
ymax=0.01,
ytick={0.0001,  0.001,   0.01},
yminorticks=false,
ylabel={$\|\vExt - \vExt^\mathrm{DEC}\|_{L^{2,2}}$},
x label style={font=\footnotesize,at={(axis description cs:0.5,0.00)}},
y label style={font=\footnotesize,at={(axis description cs:0.0,0.5)}},
yticklabel style = {font=\footnotesize},
xticklabel style = {font=\footnotesize},
axis background/.style={fill=white}
]
\addplot [color=red,solid,line width=1.0pt,mark=o,mark options={solid},forget plot]
  table[row sep=crcr]{%
32	0.00886667719576599\\
64	0.0054212661246164\\
128	0.00303478010386611\\
256	0.00162270746490437\\
512	0.000859728528543745\\
1024	0.000474652892077214\\
2048	0.000298037557779582\\
3000	0.000251450703300054\\
4096	0.000229173122578069\\
8192	0.000207806969734878\\
};
\addplot [color=black,dashed,line width=1.0pt,forget plot]
  table[row sep=crcr]{%
256	0.00390625\\
1024	0.0009765625\\
};
\addplot [color=blue,solid,line width=1.0pt,mark size=4.3pt,mark=diamond*,mark options={solid,fill=blue},forget plot]
  table[row sep=crcr]{%
3000	0.000251450703300054\\
};
\node[right, align=left, text=black]
at (axis cs:600,0.002) {\footnotesize$\sim \alpha^{-1}$};
\end{axis}
\end{tikzpicture}%
		\end{minipage}
		\begin{minipage}{0.49\textwidth}
%
%
\begin{tikzpicture}

\begin{axis}[%
width=0.7\textwidth,
height=0.5\textwidth,
at={(0,0)},
scale only axis,
separate axis lines,
every outer x axis line/.append style={black},
every x tick label/.append style={font=\color{black}},
xmode=log,
xmin=10,
xmax=10000,
xtick={   10,   100,  1000, 10000},
xminorticks=false,
xlabel={$\alpha$},
every outer y axis line/.append style={black},
every y tick label/.append style={font=\color{black}},
ymin=2.24,
ymax=2.36,
ytick={2.24,  2.3, 2.36},
ylabel={$\|\vExt^\star\cdot\surfNormal\|_{L^{2,2}} \ / \ 10^{-5}$},
x label style={font=\footnotesize,at={(axis description cs:0.5,0.0)}},
y label style={font=\footnotesize,at={(axis description cs:0.0,0.5)}},
yticklabel style = {font=\footnotesize},
xticklabel style = {font=\footnotesize},
axis background/.style={fill=white}
]
\addplot [color=red,solid,line width=1.0pt,mark=o,mark options={solid},forget plot]
  table[row sep=crcr]{%
32	2.35866759853037\\
64	2.29723977867975\\
128	2.27132037053956\\
256	2.25935139359514\\
512	2.25359669795153\\
1024	2.25078052617152\\
2048	2.24940071885282\\
3000	2.24897688376881\\
4096	2.24874354694255\\
8192	2.24847503921651\\
};
\addplot [color=blue,solid,line width=1.0pt,mark size=4.3pt,mark=diamond*,mark options={solid,fill=blue},forget plot]
  table[row sep=crcr]{%
3000	2.24897688376881\\
};
\end{axis}
\end{tikzpicture}%
		\end{minipage}
	\end{center}
 	\caption{$L^{2,2}$ norm of the error between the present velocity field $\extend{\vb}$ and the velocity field $\extend{\vb}^{\mathrm{DEC}}$ computed with DEC against the penalty parameter $\alpha$ (left) and $L^{2,2}$ norm of the normal component of the rescaled velocity field $\vExt^*=\vExt/\|\vExt\|_{L^2}$ against the penalty parameter $\alpha$ (right). The first superscript index denotes the $L^p$ norm regarding time $t$ and the second superscript index denotes the spatial $L^p$ norm. The blue diamond indicates the penalty parameter $\alpha$ used for visualization in figure \ref{fig2} and in the following examples.}
	\label{fig3}
\end{figure}

In figure \ref{fig3} (left) we compare $\vExt$ with $\extend{\vb}^{\mathrm{DEC}}$ for various $\alpha$. Thereby $\extend{\vb}^{\mathrm{DEC}}$ is the solution $\vb$ of eqs. \eqref{eq1} and \eqref{eq2} with zero normal component from \cite{Nitschkeetal_book_2017}. Again first order convergence in $\alpha$ can be obtained. In figure \ref{fig3} (right) we consider the rescaled velocity field $\vExt^* = \vExt / \| \vExt \|_{L^2}$ in order to show that the penalization of the normal component $\vExt^*\cdot\surfNormal$ is numerically satisfied. For $\alpha = const$ the same convergence properties in space and time are found as in flat geometries with periodic boundary conditions \cite{Chorin_MC_1969}.

%

\section{Results}
\label{sec4}

The Poincar\'e-Hopf theorem relates the topology of the surface to analytic properties of a vector field on it. For vector fields $\vb \in \Tangent \surf$ with only finitely many zeros (defects)
it holds that $\sum_{\xb \in \vb^{-1}(\mathbf{0})} \text{Ind}_{\xb} \vb = 2 - 2 g(\surf)$ with $\text{Ind}_\xb \vb$ the index or winding number of $\xb$ for $\vb$ and $g(\surf)$ the genus of the surface $\surf$. To highlight this relation we consider $n$-tori for $n = 1, 2, 3$ with genus $1, 2$ and $3$, respectively. Obviously, the simulation results have to fulfill the Poincar\'e-Hopf theorem in each time step, but they will also provide a realization of the theorem which depends on geometric properties and initial condition. Similar relations have already been considered for surfaces with $g(\surf) = 0$ in \cite{Reutheretal_MMS_2015,Nitschkeetal_book_2017}.

A general form of a levelset function for a $n$-torus can be written as $L(\xb) = \prod_{i=1}^nT(\xb-\mb_i) - \left( n - 1 \right)\delta$ with a constant $\delta > 0$ and the midpoints of the tori $\mb_i\in\R^3$ for $i=1,\dots,n$. In the following examples we consider the fully discrete scheme for \autoref{prob:2} and use $\text{Re} = 10$, $\tau = 0.1$, $\alpha = 3000$, $R = 1$ and $r = 0.5$. For the Gaussian curvature $\gaussianCurvature$ we use the analytic formula. The initial condition is considered to be $\vb_0 = \RotSurf \psi_0 = \surfNormal \times \GradSurf \psi_0$ with $\psi_0 = \frac{1}{2}\left(x + y + z\right)$ which ensures the incompressibility constraint. 

\begin{figure}
	\begin{center}
		\includegraphics[width=0.15\textwidth]{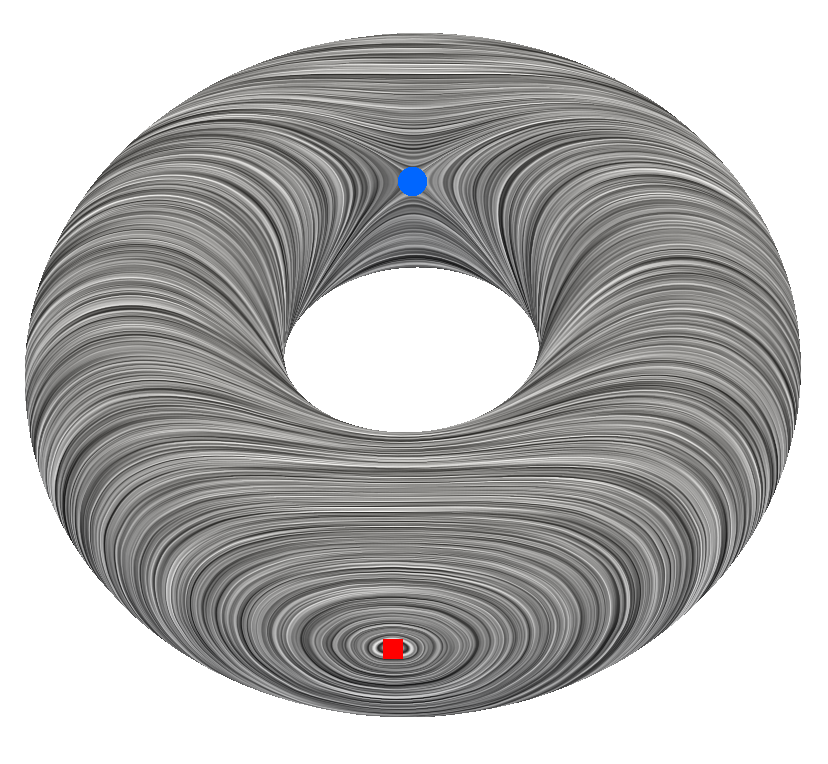}
		\includegraphics[width=0.15\textwidth]{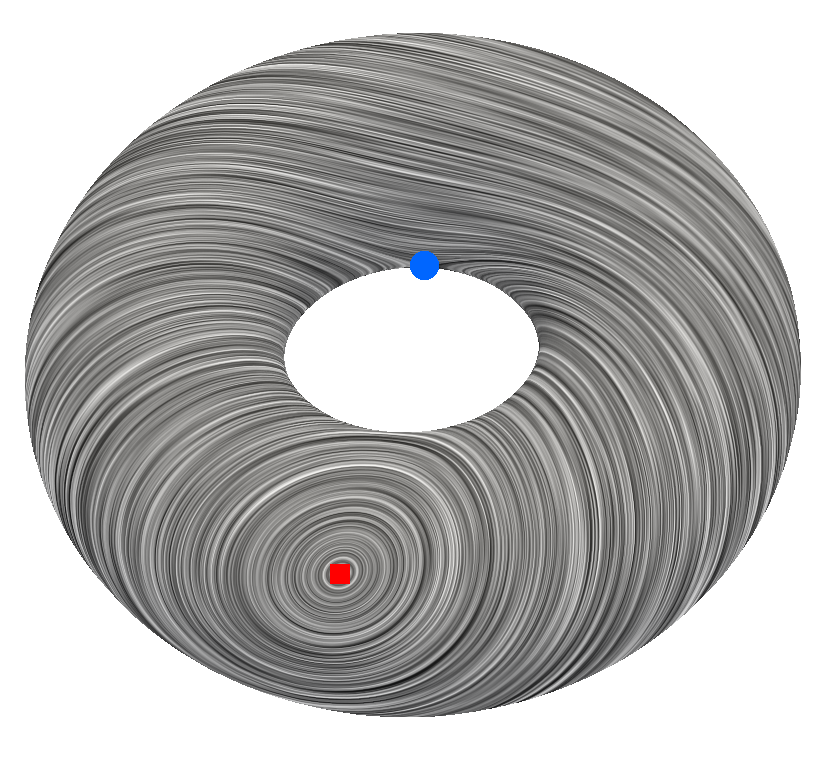}
		\includegraphics[width=0.15\textwidth]{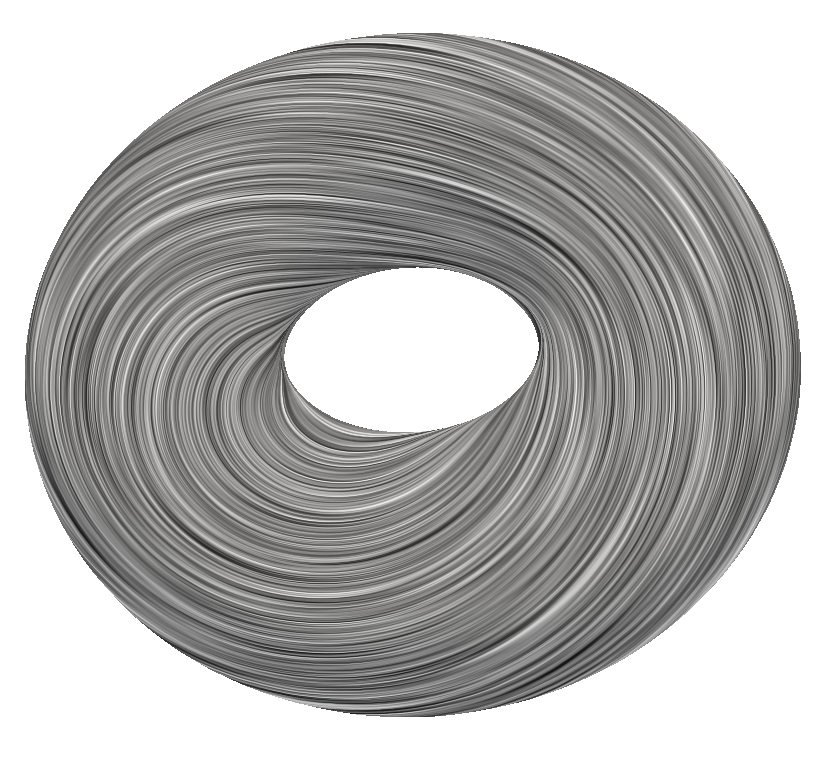}
		\includegraphics[width=0.15\textwidth]{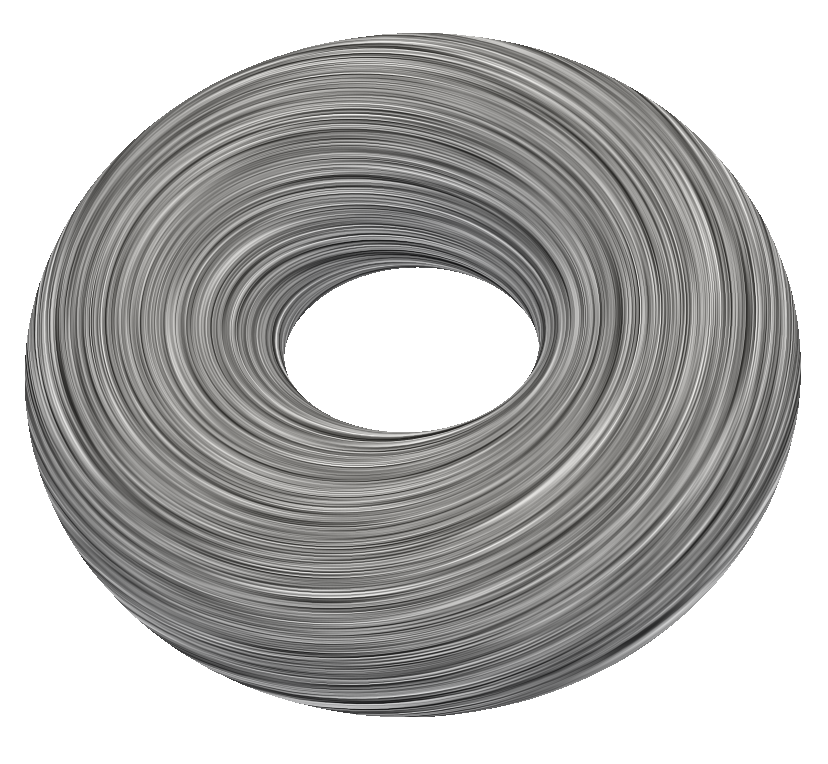}
		\includegraphics[width=0.15\textwidth]{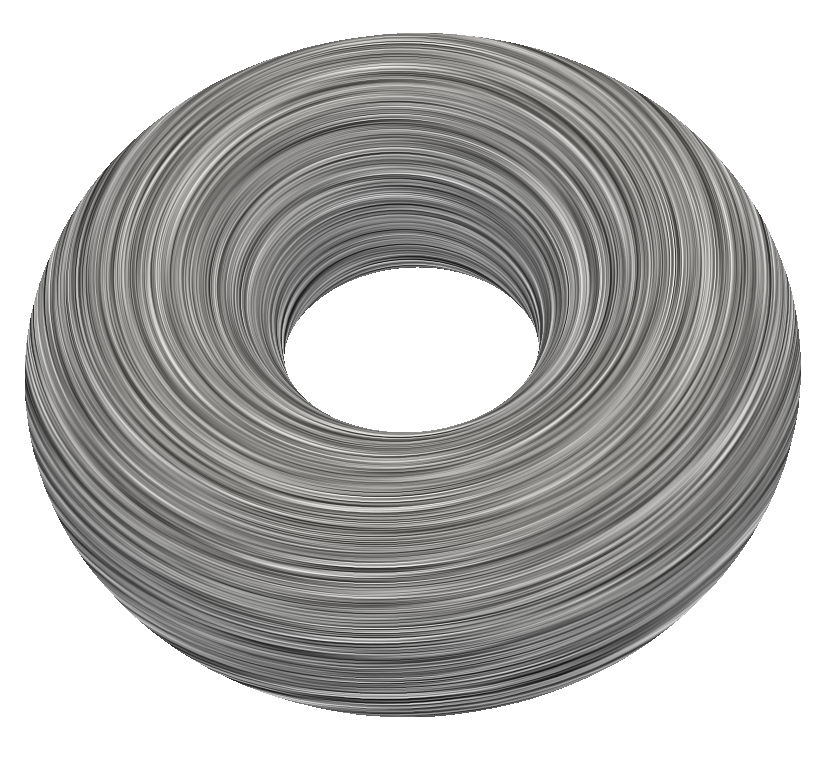}
		\includegraphics[width=0.15\textwidth]{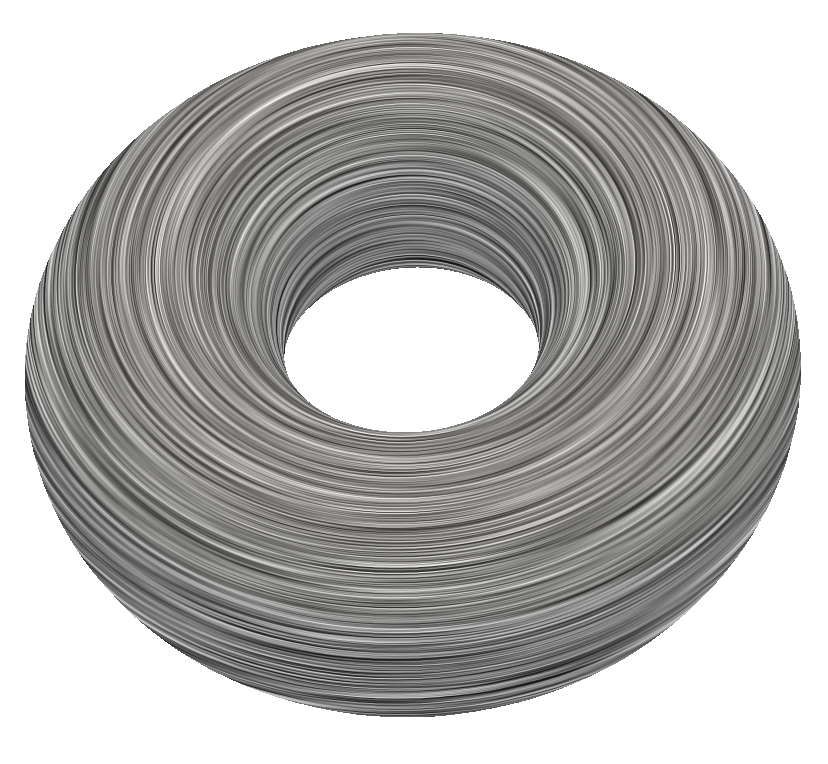}\\
		\includegraphics[width=0.26\textwidth]{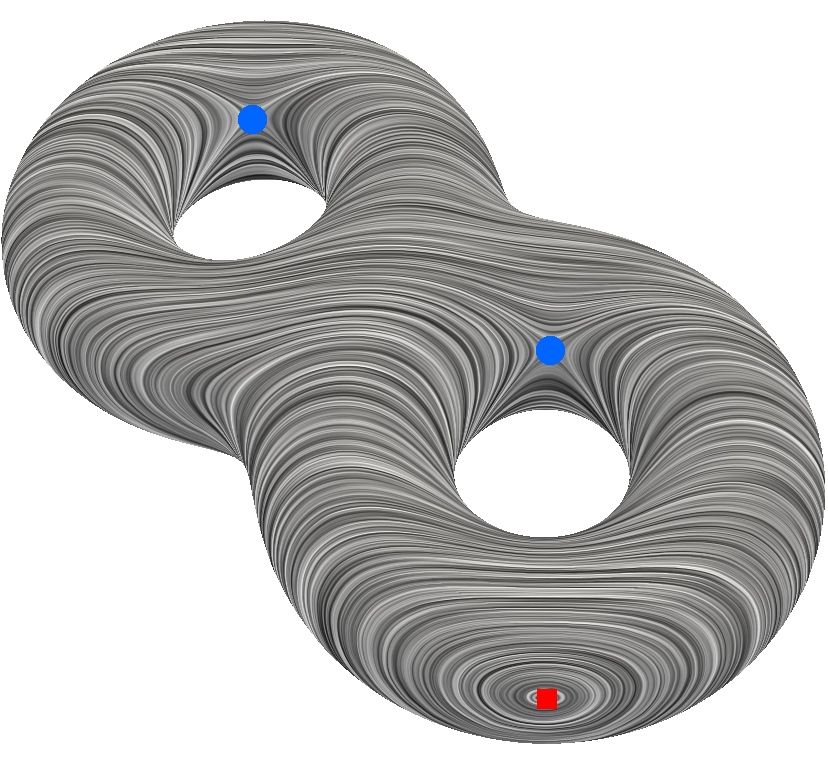}
		\includegraphics[width=0.26\textwidth]{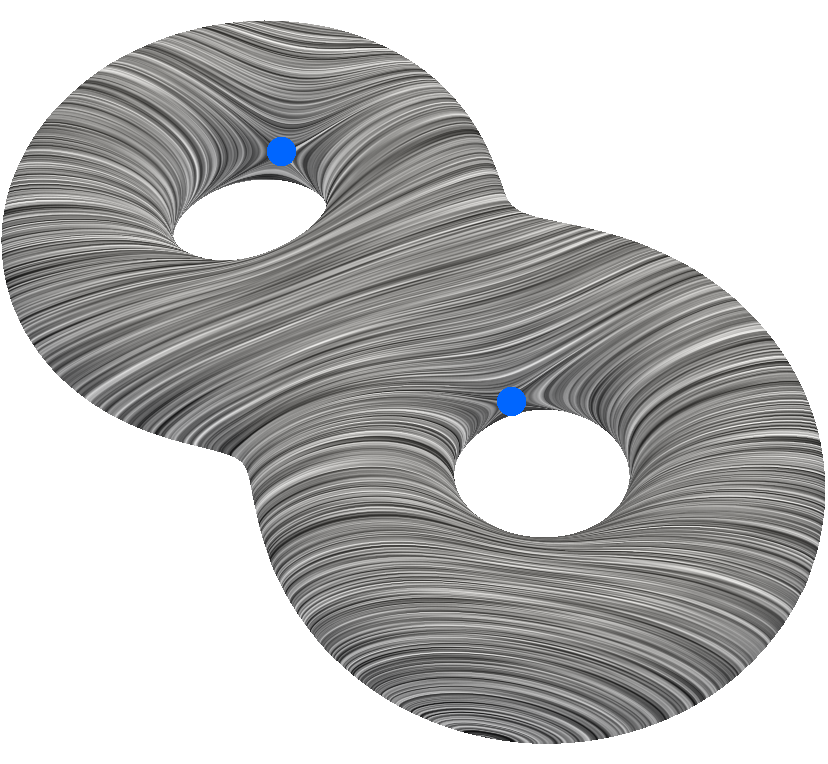}
		\includegraphics[width=0.26\textwidth]{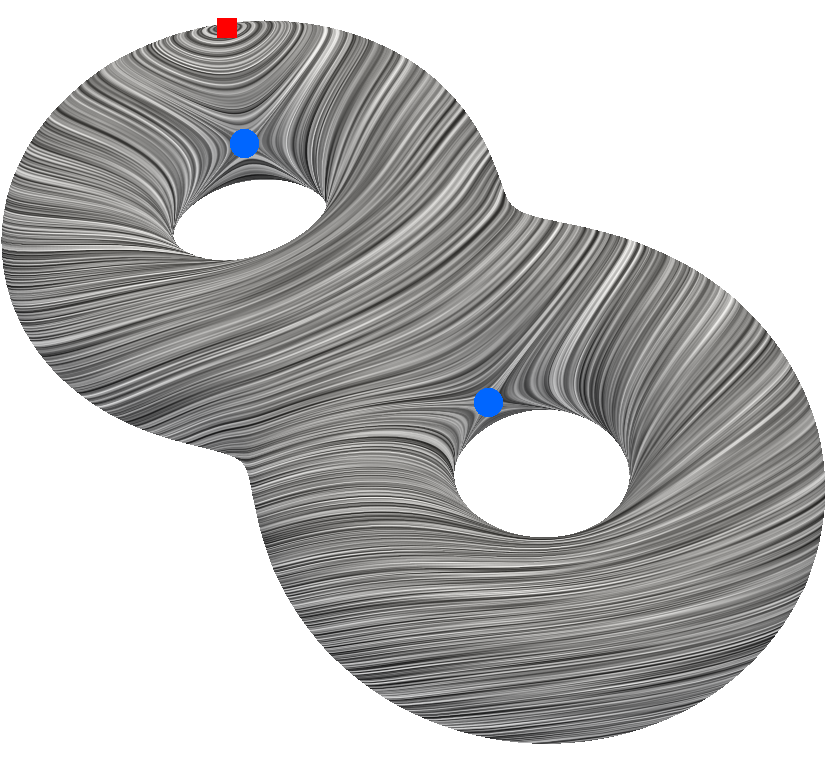}\\
		\includegraphics[width=0.26\textwidth]{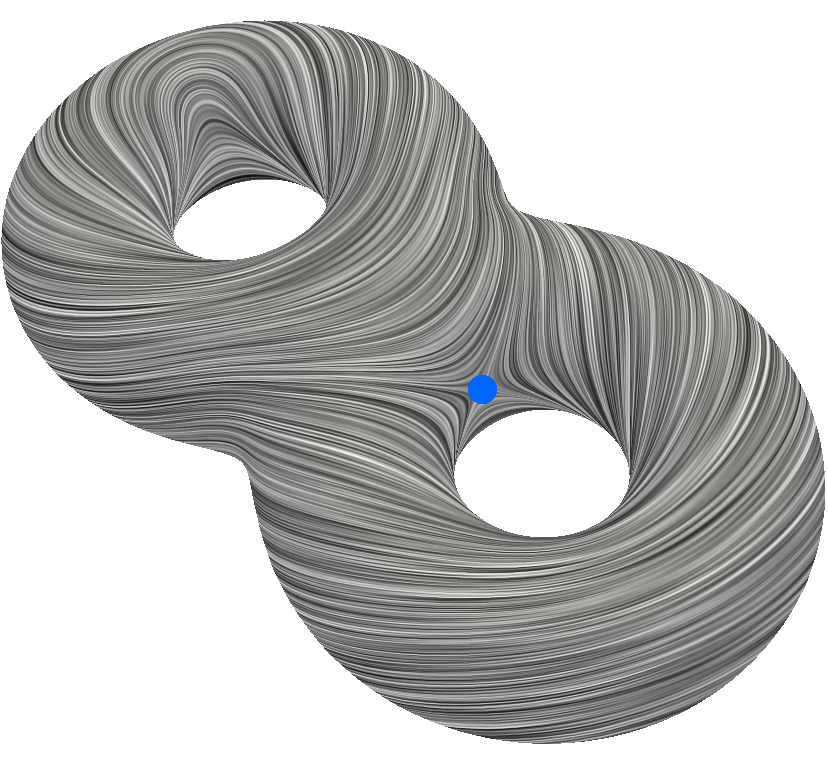}
		\includegraphics[width=0.26\textwidth]{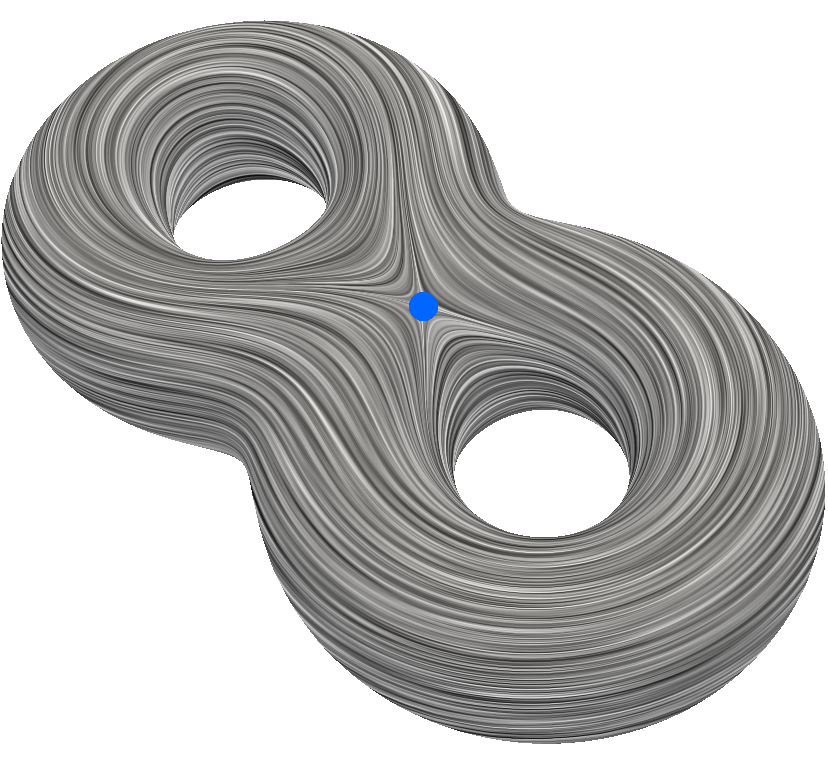}
		\includegraphics[width=0.26\textwidth]{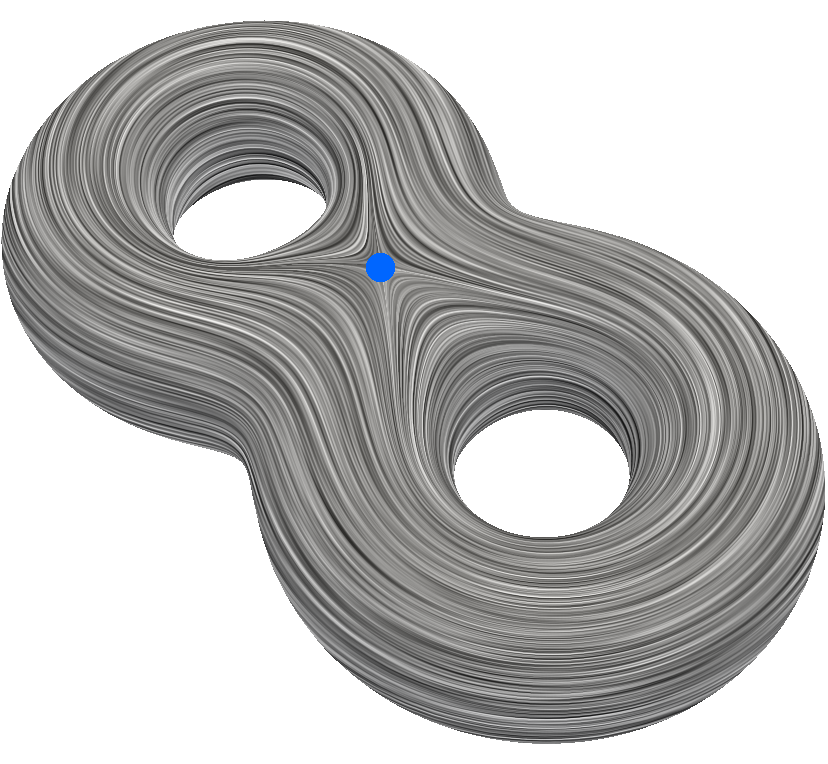}\\
		\includegraphics[width=0.26\textwidth]{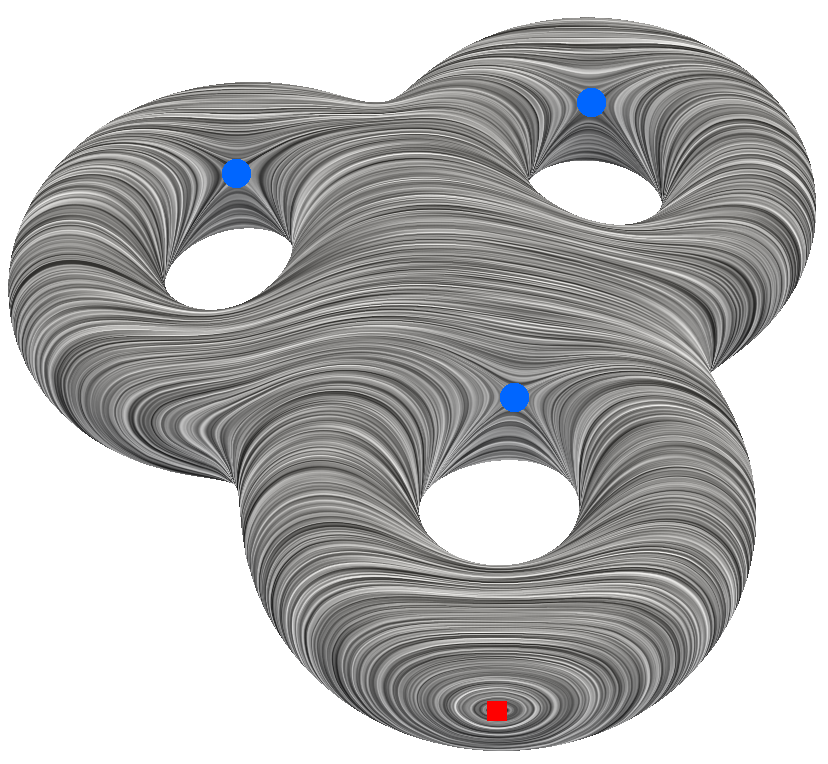}
		\includegraphics[width=0.26\textwidth]{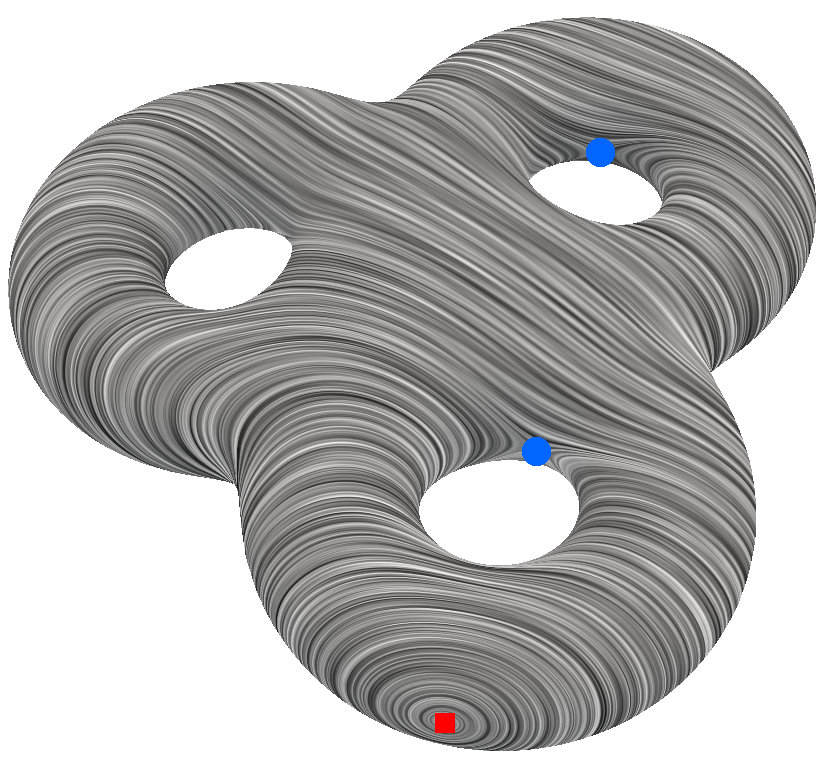}
		\includegraphics[width=0.26\textwidth]{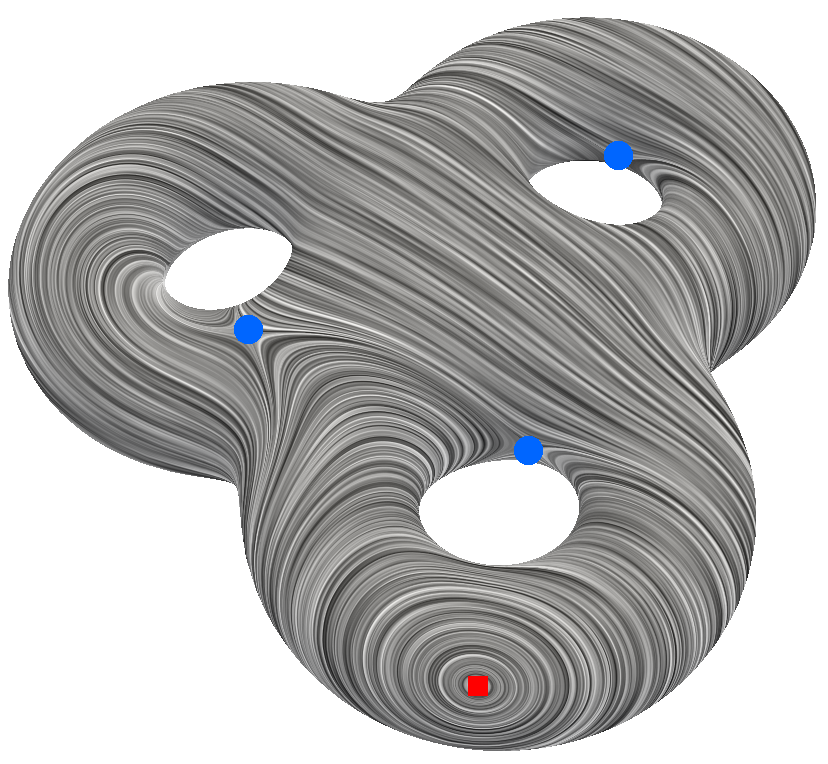}\\
		\includegraphics[width=0.26\textwidth]{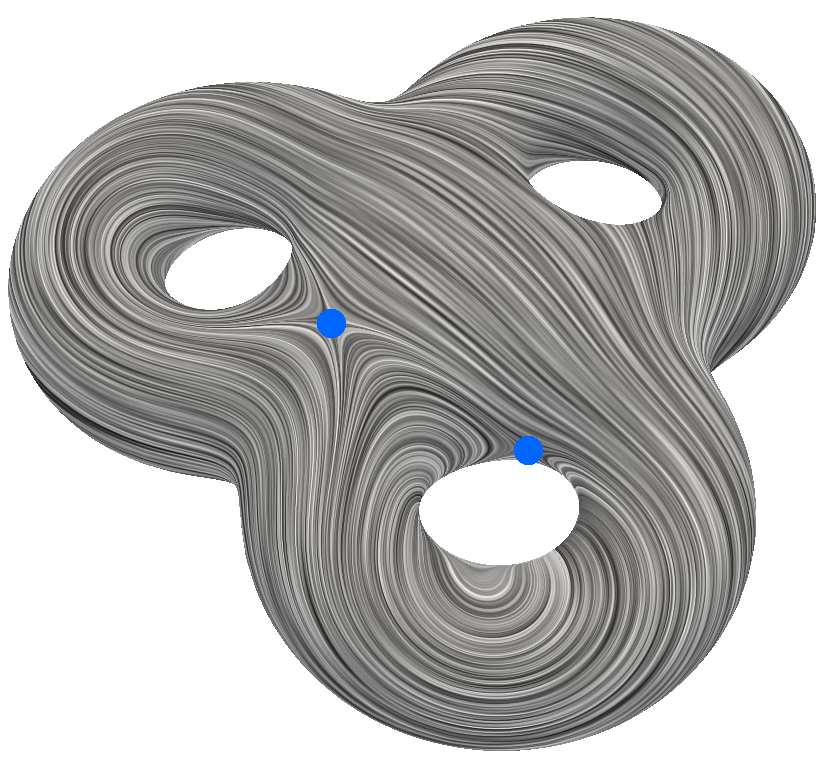}
		\includegraphics[width=0.26\textwidth]{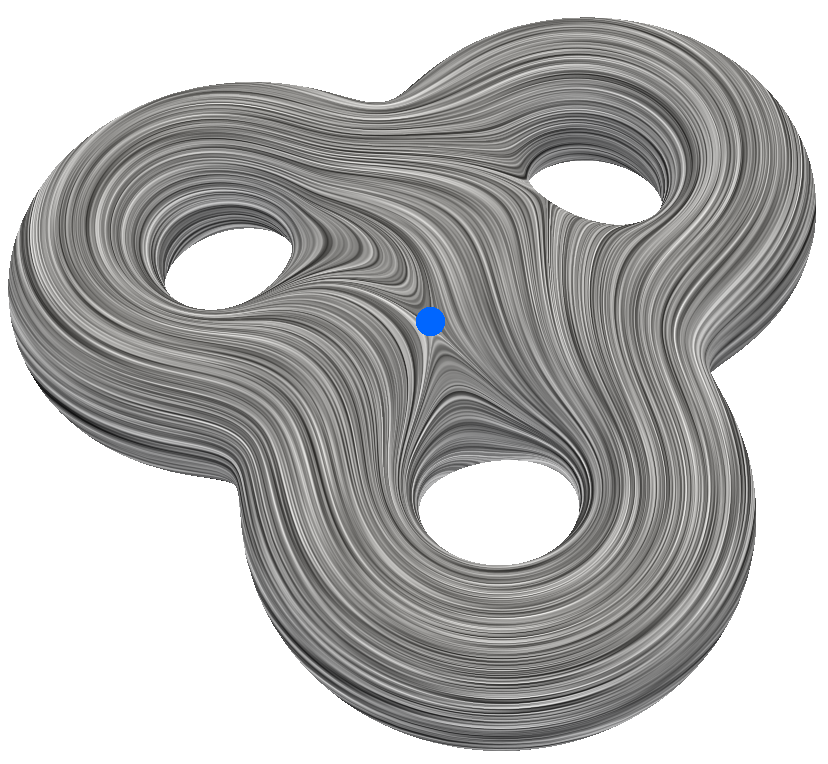}
		\includegraphics[width=0.26\textwidth]{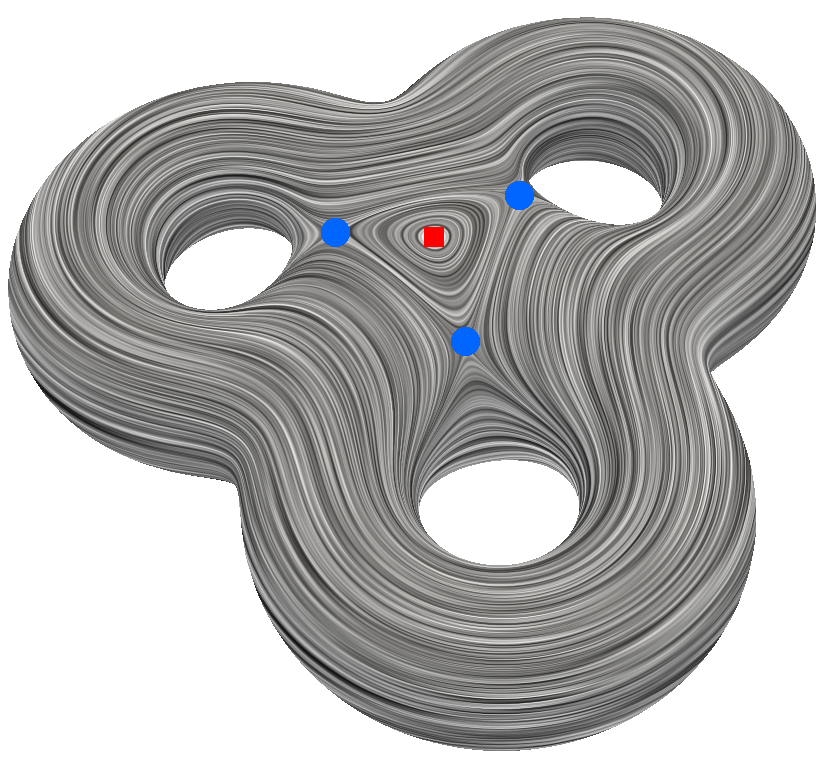}
	\end{center}
	\caption{
		Numerical solution of $\vExt = - \surfNormal \times \wExt$ for the $1$-torus (top row) at $t = 0$, $5$, $10$, $15$, $25$ and $100$ (left to right), the $2$-torus (2nd and 3rd row) at $t = 0$, $10$, $20$, $30$, $50$ and $100$ (left to right) as well as the $3$-torus (4th and 5th row) at $t = 0$, $10$, $20$, $30$, $50$ and $100$ (left to right) visualized as noise concentration field aligned to the velocity field $\vExt$. The red squares and blue circles are indicating $+1$ defects (vortices) and $-1$ defects (saddles), respectively. The full evolution for the three examples is provided in the supplementary material.
	}
	\label{fig:nTori}
\end{figure}

Figure \ref{fig:nTori} (top) shows the time evolution on the $1$-torus with $\mb_1 = \mathbf{0}$. The initial state has four defects, two vortices with $\text{Ind}_\xb \vb = +1$, indicated as red dots, and two saddles with $\text{Ind}_\xb \vb = -1$, indicated as blue dots (one vortex and one saddle are not visible). These defects annihilate during the evolution. The final state is again a Killing vector field without any defects.

For $n > 1$ the rotational symmetry is broken and Killing vector fields are no longer possible. We thus expect dissipation of the kinetic energy and convergence to $\vb = \mathbf{0}$ for any initial condition.
Figure \ref{fig:nTori} (middle) shows the time evolution on a $2$-torus where we have used the midpoints $\mb_1 = (-1.2, 0, 0)$ and $\mb_2 = - \mb_1$ as well as $\delta = 1$. The initial state has two vortices and four saddles and thus $\sum_{\xb \in \vb^{-1}(\mathbf{0})} \text{Ind}_{\xb} \vb = -2$. Two vortex-saddle pairs annihilate each other and the final defect configuration consists of two saddles located at the center of the $2$-torus (one is not visible). The velocity field decays towards $\vb = \mathbf{0}$.
Figure \ref{fig:nTori} (bottom) shows the time evolution on a $3$-torus with midpoints $\mb_1 = (-1.2, -0.75, 0)$, $\mb_2 = (1.2, -0.75, 0)$ and $\mb_3 = (0, 1.33, 0)$ as well as $\delta = 10$. Initially we have three vortices and seven saddles and thus $\sum_{\xb \in \vb^{-1}(\mathbf{0})} \text{Ind}_{\xb} \vb = -4$, which is also fulfilled for the final defect configuration with two vortices and six saddles at the center of the $3$-torus (one vortex and three saddles are not visible). Again the velocity field decays towards $\vb = \mathbf{0}$.

\begin{figure}
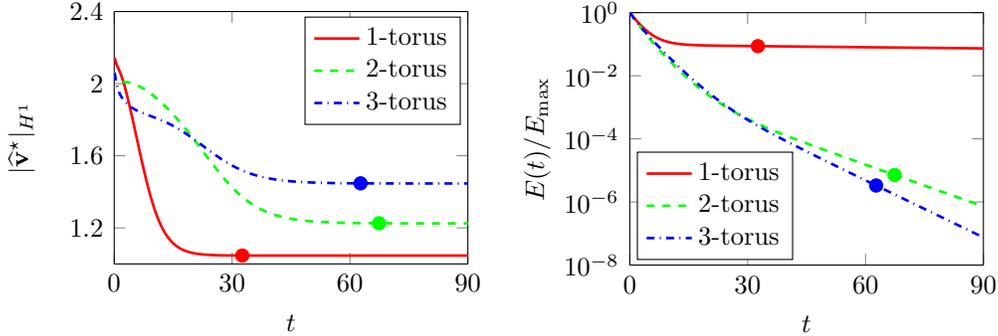

	\begin{minipage}{0.49\textwidth}
		\begin{center}
			\input{h1SemiNorm.tex}
		\end{center}
	\end{minipage}
	\begin{minipage}{0.49\textwidth}
		\begin{center}
			\input{nToriKineticEnergy.tex}
		\end{center}
	\end{minipage}
	\caption{$H^1$ semi-norm of the rescaled velocity field $\vExt^*=\vExt/\|\vExt\|_{L^2}$ against time $t$ (left) and normalized kinetic energy $E/E_{\mathrm{max}}$ against time $t$ (right), where $E_{\mathrm{max}}$ is the maximum value of the kinetic energy $E$ over time. The colored dots are indicating the time points at which the defects reach their final position and only viscous dissipation takes place or a Killing vector field is formed. We identify these points if the decay rate of the $H^1$ semi-norm of the rescaled velocity field $\vExt^*$ reaches $0.001\%$ of its maximum value over time.}
	\label{fig:toriComparison}
\end{figure}

To show the differences in the evolution on the $n$-tori before and after the final defect configuration is reached we consider the $H^1$ semi-norm of the rescaled velocity field $\vExt^* = \vExt / \| \vExt \|_{L^2}$. If the defects do not move this quantity is constant. Figure \ref{fig:toriComparison} shows the evolution over time together with the decay of the kinetic energy $E = \frac{1}{2} \int_\surf \|\vExt\|^2 \; d \surf$. 

These results clearly show the strong interplay between topology, geometric properties and defect positions. 
	
\section{Conclusions}
\label{sec5}

We have proposed a discretization approach for the incompressible surface Navier-Stokes equation on surfaces with arbitrary genus $g(\surf)$. The approach only requires standard ingredients which most finite element implementations can offer. It is based on a reformulation of the equation in Cartesian coordinates of the embedding $\R^3$, penalization of the normal component, a Chorin projection method and discretization is space by globally continuous, piecewise linear Lagrange surface finite elements for each component. A further rotation of the velocity field leads to a drastic reduction of the complexity of the equation and the required computing time. The fully discrete scheme is described in detail and its accuracy validated against a DEC solution on a $1$-torus, which was considered in \cite{Nitschkeetal_book_2017}. The interesting interplay between the topology of the surface, its geometric properties and defects in the flow field are shown on $n$-tori for $n = 1,2,3$. 

Even if the formulation of the incompressible surface Navier-Stokes equation is relatively old \cite{Scriven_CES_1960,EbinMarsden_AM_1970,MitreaTaylor_MA_2001}, numerical treatments on general surfaces are very rare. We are only aware of the DEC approach in \cite{Nitschkeetal_book_2017} and therefore expect the proposed approach to initiate a broader use and advances in the mentioned applications in Section \ref{sec1}. We further expect it to be the basis for further developments, e.g. coupling of the surface flow with bulk flow in two-phase flow problems, as, e.g. considered in \cite{Reutheretal_JCP_2016} using a vorticity-stream function approach or in \cite{Barrettetal_NM_2016} within an alternative formulation based on the bulk velocity and projection operators. Another extension considers evolving surfaces. With a prescribed normal velocity this has already been considered in \cite{Reutheretal_MMS_2015}, again using a vorticity-stream function approach. The corresponding equations are derived in \cite{kobaetal_QAM_2017} using a global variational approach and in \cite{Miura_arXiv_2017} as a thin-film limit. A mathematical derivation of the evolution equation for the normal component is still controversial. The derivation in \cite{ArroyoDeSimone_PRE_2009} is based on local conservation of mass and linear momentum in tangential and normal direction, while the derivation in \cite{Jankuhnetal_arXiv_2017} is based on local conservation of mass and total linear momentum. The resulting equations differ. However, in the special case of a stationary surface, all these models coincide with the incompressible surface Navier-Stokes equation in eqs. \eqref{eq1} and \eqref{eq2}.

{\bf Acknowledgements:} This work is partially supported by the German Research Foundation through grant Vo899/11. We further acknowledge computing resources provided at JSC under grant HDR06 and at ZIH/TU Dresden.

\bibliographystyle{elsarticle-harv}
\bibliography{biball}

\end{document}